\input amstex
\documentstyle{amsppt}
\magnification=\magstep1
\hsize=5.2in
\vsize=6.8in

\centerline{\bf ON RELATIVE PROPERTY (T) AND HAAGERUP'S PROPERTY}

\vskip 0.05in

\vskip 0.2in
\centerline {\rm by}
\vskip 0.05in
\centerline {\rm IONUT CHIFAN and ADRIAN IOANA\footnote{The second author was supported by a Clay Research Fellowship}}

\address Math Dept, UCLA, LA, CA 90095-155505
\endaddress
\email ichifan\@math.ucla.edu\endemail
\address Math Dept  UCLA, LA, CA 90095-155505\endaddress
\email adiioana\@math.ucla.edu \endemail
\topmatter
\abstract 
 We consider the following three properties for countable discrete groups $\Gamma$: (1) $\Gamma$ has an infinite subgroup with relative property (T), (2) the group von Neumann algebra $L\Gamma$ has a diffuse von Neumann subalgebra with relative property (T) and (3) $\Gamma$ does not have Haagerup's property. It is clear that (1) $\Longrightarrow$ (2) $\Longrightarrow$ (3). We prove that both of the converses are false.
\endabstract

\endtopmatter

\document

\head \S {0. Introduction.} \endhead
\vskip 0.1in

In this paper, we investigate the relationship between {\it relative property (T)} and {\it Haagerup's property} in the context of countable groups and finite von Neumann algebras. An inclusion $(\Gamma_0\subset\Gamma)$ of countable discrete groups has {\it relative property (T)} of Kazhdan-Margulis if any unitary representation of $\Gamma$ which has almost invariant vectors necessarily has a non-zero $\Gamma_0$-invariant vector. The classical examples here are $(\Bbb Z^2\subset\Bbb Z^2\rtimes$SL$_2(\Bbb Z))$ and (SL$_n(\Bbb Z)\subset$ SL$_n(\Bbb Z)$), for $n\geq 3$ ([Ka67],[Ma82]).

The presence of  relative property (T)  subgroups is an obstruction to Haagerup's property.  
A countable group $\Gamma$ is {\it Haagerup} if it admits a $c_0$ (or {\it mixing}) unitary representation which has almost invariant vectors. This class includes amenable groups, free groups and is closed under free products and, more surprisingly, wreath products ([CSV09]).

It is clear from the definitions that a group with  Haagerup's property cannot have relative property (T) with respect to any infinite subgroup. In  [CCJJV01, Section 7.1] the authors asked whether the converse holds true, i.e. if having an infinite subgroup with relative property (T) is the {\it only} obstruction to Haagerup's property. This question was answered negatively by Y. de Cornulier ([Co06ab]). He showed that there are certain  groups, e.g. $\Gamma_{\alpha}=\Bbb Z[\root 3\of{\alpha}]^3\rtimes SO_3(\Bbb Z[\root 3\of{\alpha}])$, for $\alpha$ not a cube, which do  not have relative property (T) with respect to any infinite subgroup but yet have  relative property (T) with respect to some infinite {\it subset} (see Definition 1.1). The latter property still guarantees the failure of Haagerup's property.

The notion of {\it relative property (T)} (or {\it rigidity}) for inclusions of finite von Neumann algebras was introduced by S. Popa in [Po06a] (see Definition 1.4). Since then it has found many striking applications to  von Neumann algebras theory and orbit equivalence ergodic theory (see the surveys [Po07] and [Fu09]). Thus, it was used as a key ingredient in Popa's solution to the long standing problem of finding II$_1$ factors with trivial fundamental group ([Po06a]).  Examples of rigid inclusions of von Neumann algebras are provided by inclusions of groups. Precisely, the inclusion $(\Gamma_0\subset \Gamma)$ of two countable groups has relative property (T) if and only if the inclusion of group von Neumann algebras $(L\Gamma_0\subset L\Gamma)$ has it ([Po06a]). Similarly, there is a notion of Haagerup's property for finite von Neumann algebras which generalizes the corresponding notion for groups ([CJ85],[Ch83]). 

Relative property (T) and Haagerup's property are also incompatible in the framework of von Neumann algebras: if a finite von Neumann algebra $N$ has  Haagerup's property then it does not have any diffuse (i.e. non-atomic) relatively rigid subalgebra. The main goal of this paper is to show that converse of this statement is false even in the case of von Neumann algebras $N$ arising from countable groups. Thus, we are interested in finding a countable group $\Gamma$ such that its group von Neumann algebra $L\Gamma$  neither has Haagerup's property nor any diffuse rigid subalgebra.
The natural candidates are de Cornulier's examples, $\Gamma_{\alpha}$, since one might hope that the absence of relative property (T) subgroups is inherited by the group von Neumann algebra.
However, we prove that this is not the case: 

\proclaim{0.1 Theorem}(Corollary 2.2) Let $\alpha\in\Bbb N\setminus\{\beta^3|\beta\in\Bbb N\}$ and let $\Gamma_{\alpha}=\Bbb Z[\root 3\of{\alpha}]^3\rtimes SO_3(\Bbb Z[\root 3\of{\alpha}])$.  Then there exists a diffuse von Neumann subalgebra $B$ of $N=L(\Bbb Z[\root 3\of{\alpha}]^3)$ such that the inclusion $(B\subset L(\Gamma_{\alpha}))$ is rigid and $B'\cap N=L(\Bbb Z[\root 3\of{\alpha}]^3)$.
\endproclaim

In [Co06b], a concept of {\it resolutions} was introduced in order to quantify the transfer of property (T) from a locally compact group to its lattices. In particular, they can be used to locate the relative property (T) subsets of  lattices $\Gamma$ in Lie groups, e.g. $\Gamma=\Gamma_{\alpha}$.  The idea behind the proof of Theorem 0.1 is that by combining resolutions with techniques from [Io09] we can also detect certain rigid subalgebras of $L\Gamma$. 

Next, we consider another class of groups. If $A$ and $\Gamma$ are two countable groups and $X$ is a countable $\Gamma$-set then the {\it generalized wreath product} group $A\wr_{X}\Gamma$ is defined as $A^{X}\rtimes\Gamma$. If $X=\Gamma$, together with the left multiplication $\Gamma$-action, then we recover the standard wreath product $A\wr\Gamma$. Following results from [Po06ab] and [Io07] we know that if $A$ and $\Gamma$ have Haagerup's property, then the von Neumann algebra $L(A\wr_{X}\Gamma)$ does not have a diffuse rigid von Neumann subalgebra,  regardless of the set $X$. Thus, in order to get examples of groups with the desired properties, we just need to find a suitable set $X$ for which $A\wr_{X}\Gamma$ is not Haagerup.

\proclaim {0.2 Theorem} (Corollary 3.4) Let $A$ be a non-trivial countable Haagerup group. Let $\Gamma$ be a countable Haagerup together with a quotient group $\Gamma_0$. Assume that $\Gamma_0$ is not Haagerup and endow it with the left multiplication action of $\Gamma$. 
Then $A\wr_{\Gamma_0}\Gamma$ is not Haagerup. 
Thus, the group von Neumann algebra $N=L(A\wr_{\Gamma_0}\Gamma)$ does not have Haagerup property and does not admit any diffuse von Neumann subalgebra $B$ such that the inclusion $(B\subset  N)$ is rigid.
\endproclaim

The proof of Theorem 0.2 is based on a general result:  a semidirect product $A\rtimes\Gamma$ is not Haagerup whenever $A$ is  abelian  and $\Gamma$ acts on $A$ through a non-Haagerup quotient group, $\Gamma_0$.

Recently, de Cornulier, Stalder and Valette proved that the class of Haagerup groups is closed under standard wreath products ([CSV09]). Moreover, they showed that if $A$ and $\Gamma$ are Haagerup groups then the generalized wreath product $A\wr_{X}\Gamma$ is Haagerup for certain sets $X$ and conjectured  that this is the case for any $X$. Theorem 0.2 provides in particular a counterexample to their conjecture. For example, if $\Gamma$ is a free group and $\Gamma_{0}$ is a property (T) quotient of $\Gamma$ , then $\Bbb Z\wr_{\Gamma_0}\Gamma$ is not Haagerup. 
\vskip 0.1in
\noindent
{\it Acknowledgment}. We are grateful to Professors Sorin Popa  and Yehuda Shalom for useful discussions and encouragement.
\vskip 0.2in
\head \S {1. Preliminaries}.\endhead
\vskip 0.1in
We start this section by reviewing the notion of relative property (T) for groups. Then we explain de Cornulier's examples of groups which are not 
Haagerup but do not have any infinite subgroup with relative property (T). Finally, we recall Popa's notion of rigidity for inclusions of  von Neumann algebras.

A continuous unitary representation $\pi$ of a locally compact group $G$ on a Hilbert space $\Cal H$  has {\it almost invariant vectors} if for all $\varepsilon>0$ and any compact set $F\subset G$ we can find a unit vector $\xi\in\Cal H$ such that $||\pi(g)(\xi)-\xi||\leq \varepsilon,$ for all $g\in F$. If $H$ is a closed subgroup of $G$ then the inclusion $(H\subset G)$ has {\it relative property }(T) of Kazhdan-Margulis if any  unitary representation of $G$ which has almost invariant vectors must have a non-zero $H$-invariant vector ([K67],[Ma82]).
Recall also that $G$ has {\it Haagerup's property} (in short, {\it is Haagerup}) if it admits a $c_0$  unitary representation $\pi:G\rightarrow\Cal U(\Cal H)$ which has almost invariant vectors. Being $c_0$ means that for every $\xi,\eta\in\Cal H$ we have that $\lim_{g\rightarrow\infty}\langle\pi(g)(\xi),\eta\rangle=0$.

If a countable, discrete group $\Gamma$ is Haagerup then it does not have relative property (T) with respect to any infinite subgroup. de Cornulier  proved that the converse is false ([Co06ab]).  For example, he showed that if $\alpha\in\Bbb N\setminus\{\beta^3|\beta\in\Bbb N\}$, then   $\Gamma_{\alpha}=\Bbb Z[\root 3\of{\alpha}]^3\rtimes SO_3(\Bbb Z[\root 3\of{\alpha}])$
neither has Haagerup property nor admits an infinite subgroup with relative property (T). 
To quickly see that $\Gamma_{\alpha}$ is not Haagerup just notice that it is measure equivalent (see [Fu09] for the definition) to the group $\Lambda=\Bbb Z^3\times (\Bbb Z[i]^3\rtimes SO_3(\Bbb Z[i]))$ which has an infinite subgroup with relative property (T)  (i.e. $\Bbb Z[i]^3$). Indeed, both $\Gamma_{\alpha}$ and $\Lambda$ are lattices in $G=(\Bbb R^3\rtimes SO_3(\Bbb R))\times (\Bbb C^3\rtimes SO_3(\Bbb C))$ ([Ma91],[Wi08]). This example shows that having an infinite subgroup with relative property (T) is not a measure equivalence invariant. 
To better explain the failure of Haagerup's property for $\Gamma_{\alpha}$,  the following two notions were introduced in [Co06b]:

\vskip 0.1in

\noindent 
{\bf 1.1 Definitions [Co06b]} (a)  Let $p:\Gamma\rightarrow G$ be a morphism between two locally compact groups with dense image.  We say that $p$ is a {\it resolution} if for any unitary representation $\pi$ of $\Gamma$ which has almost invariant vectors, there exists a subrepresentation $\sigma$ of $\pi$ of the form $\sigma=\tilde{\sigma}\circ p$, where $\tilde\sigma$ is a unitary representation of $G$ which has almost invariant vectors.
\vskip 0.02in
\noindent 
(b) Given a subset $X$ of a locally compact group $G$, we say that the inclusion $(X\subset G)$ has {\it relative property (T)} if for any unitary representation $\pi:G\rightarrow\Cal U(\Cal H)$  which has almost invariant vectors and  any $\varepsilon>0$, we can find a unit vector $\xi\in\Cal H$ such that  $||\pi(g)(\xi)-\xi||\leq \varepsilon,$ for all $g\in X$.
\vskip 0.05in
Resolutions are  useful to encode the transfer of relative property (T) from a group to its lattices. To see this, assume that $\Gamma$ is a lattice in a locally compact group $G$, let $H\subset G$ be a normal, closed subgroup such that the inclusion $(H\subset G)$ has relative property (T) and let $p:G\rightarrow G/H$ be the projection. Under these assumptions,  [Co06b, Theorem 4.3.1] asserts that  the morphism $p_{|\Gamma}:\Gamma\rightarrow \overline{p(\Gamma)}$ is a resolution. In the case when $G$ has property (T) and $H=G$, this is just saying any lattice $\Gamma$ of $G$ has property (T), thus recovering Kazhdan's classical result ([K67]).

Now, if $G=(\Bbb R^3\rtimes SO_3(\Bbb R))\times (\Bbb C^3\rtimes SO_3(\Bbb C))$ and $H=\Bbb C^3$, then  the inclusion $(H\subset G)$ has relative property (T) (see [Co06b, 3.3.1]). 
By applying  the above theorem to this situation the following was deduced in [Co06b, the proof of 4.6.3]:
\proclaim {1.2 Corollary [Co06b]} The inclusion $\Gamma_{\alpha}\hookrightarrow \Bbb R^3\rtimes SO_3(\Bbb Z[\root 3\of{\alpha}])$ is a resolution. Thus, if $\Cal B$ is the unit ball of $\Bbb R^3$ and $X=\Bbb Z[\root 3\of{\alpha}]^3\cap \Cal B$, then $(X\subset \Gamma_{\alpha})$ has relative property (T), for every $\alpha$. In particular, $\Gamma_{\alpha}$ is not Haagerup.
\endproclaim
\vskip 0.05in
Notice moreover that $X$ is a normal subset of $\Gamma_{\alpha}$.
In relation to this, let us note that results from [Co06b] imply that {\it any} lattice in a connected Lie group either has Haagerup's property or admits an infinite, ``almost normal" subset with relative property (T):

\proclaim {1.3 Corollary} Let $G$ be a connected Lie group which does not have Haagerup's property. Let $\Gamma$ be a lattice in $G$.  Then there exists an infinite set $X\subset\Gamma$ such that the inclusion $(X\subset\Gamma)$ has relative property (T) and $\gamma X\gamma^{-1}\cap X$ is infinite, for all $\gamma\in\Gamma$.
\endproclaim
{\it Proof.} Since $G$ is not Haagerup, by [Co06b, 3.3.1] and [CCJJV01, Chapter 4] we get that it has a non-trivial, normal, closed subgroup $H$ such that  the inclusion $(H\subset G)$ has relative property (T) and $G/H$ has Haagerup's property. Let $p:G\rightarrow G/H$ denote the projection and set $Q=\overline{p(\Gamma)}$. By [Co06b, 4.3.1]   the morphism $p_{|\Gamma}:\Gamma\rightarrow Q$ is a resolution. 
Thus, since the inclusion $(\{1\}\subset Q)$ has relative property (T) by [Co06b, 4.2.6] we deduce that the inclusion $((\Gamma\cap H)\subset \Gamma)$ also does. 

Since $\Gamma\cap H$ is a normal subgroup of $\Gamma$ we can hereafter assume that it is finite (otherwise, we can take $X=\Gamma\cap H$). Under this assumption, we claim that $p(\Gamma)$ is a non-discrete subgroup of $G/H$. Indeed, if $p(\Gamma)$ is discrete then it must have Haagerup's property, as $G/H$ has it. But $p(\Gamma)$ is isomorphic to $\Gamma/(\Gamma\cap H)$ and since $\Gamma\cap H$ is finite, we would get that $\Gamma$ is Haagerup, a contradiction.

Next, let $V$ be a neighborhood of $1\in Q$ with compact closure and define $X=p^{-1}(V)\cap \Gamma$.  
Since the inclusion $(V\subset Q)$ has relative property (T) by [Co06b, 4.2.6] we deduce that the inclusion $(X\subset \Gamma)$ has relative property (T). To check the normality assertion, fix $\gamma\in\Gamma$ and denote $W=p(\gamma)Vp(\gamma)^{-1}\cap V$. Then $\gamma X\gamma^{-1}\cap X=\{x\in\Gamma|p(x)\in W\}.$ Since $p(\Gamma)\subset G/H$ is non-discrete and $W$ is a neighborhood of $1\in Q$, the latter set is infinite.
\hfill$\blacksquare$ 
\vskip 0.05in

To remind Popa's notion of rigidity for von Neumann algebras, let $M$ be a separable finite von Neumann algebra with a faithful, normal trace $\tau:M\rightarrow\Bbb C$ and let $B\subset M$ be a  von Neumann subalgebra. A Hilbert space $\Cal H$ is called a  {\it Hilbert $M$-bimodule} if it admits commuting left and right Hilbert $M$-module structures. A vector $\xi\in\Cal H$ is called {\it tracial} if $\langle x\xi,\xi\rangle=\langle\xi x,\xi\rangle=\tau(x)$, for all $x\in M$, and {\it $B$-central} if $b\xi=\xi b,$ for all $b\in B$.  A Hilbert $M$-bimodule $\Cal H$ together  with a unit vector $\xi\in\Cal H$ is called a {\it pointed Hilbert $M$-bimodule} and is denoted $(\Cal H,\xi)$.
\vskip 0.05in
\noindent
{\bf 1.4 Definition [Po06a]} The inclusion $(B\subset M)$ is {\it rigid} (or has {\it relative property (T)}) if for every $\varepsilon>0$ there exists $F\subset M$ finite and $\delta>0$ such that whenever $(\Cal H,\xi)$ is a pointed   Hilbert $M$-bimodule with $\xi$ a tracial vector verifying $||x\xi-\xi x||\leq\delta,$ for all $x\in F$, there exists a $B$-central vector $\eta\in\Cal H$ with $||\eta-\xi||\leq\varepsilon$.
\vskip 0.05in
The notion of rigidity for inclusions of von Neumann algebras is analogous to and generalizes the notion of relative property (T) for groups.
More precisely, given two countable groups $\Gamma_0\subset\Gamma$, the inclusion $(L(\Gamma_0)\subset L(\Gamma))$ of their group von Neumann algebras is rigid if and only if the inclusion $(\Gamma_0\subset\Gamma)$ has relative property (T) ([Po06a, Proposition 5.1.]). Now, if $\Gamma$ is a non-amenable subgroup of SL$_2(\Bbb Z)$ acting on $\Bbb Z^2$ by matrix multiplication, then the inclusion $(\Bbb Z^2\subset\Bbb Z^2\rtimes\Gamma)$ has relative property (T) ([Bu91, section 5]) and therefore the inclusion $(L(\Bbb Z^2)\subset L(\Bbb Z^2\rtimes\Gamma))$ is rigid.

The first examples of rigid inclusions of von Neumann algebras which do not rely relative property (T) for some pair of groups have been recently exhibited in [Io09]. Thus, it is shown that for {\it any} non-amenable subfactor $N$ of $L(\Bbb Z^2\rtimes$SL$_2(\Bbb Z))$ which contains $L(\Bbb Z^2)$, we have that the inclusion $(L(\Bbb Z^2)\subset N)$ is rigid ([Io09, Theorem 3.1]).

\vskip 0.2in

\head  \S {2. Rigid subalgebras from resolutions.}\endhead
\vskip 0.1in
The main goal of this section is to show that, in certain situations, resolutions can be used to construct rigid subalgebras of von Neumann algebras (Theorem 2.1). Thus, we employ the resolution provided by Corollary 1.2 to deduce that the group von Neumann algebra $L(\Gamma_{\alpha})$ has a diffuse rigid subalgebra (Corollary 2.2). This result should be contrasted with the fact that $\Gamma_{\alpha}$ has no infinite subgroup with relative property (T).

\proclaim {2.1 Theorem} Let $\Gamma$ be a countable subgroup of SO$_n(\Bbb R)$, for some $n\geq 3$, and consider the natural action of $\Gamma$ 
on $H=\Bbb R^n$.

$\bullet$ Assume that  $A\simeq\Bbb Z^m$ ($m\geq n+1$) is a $\Gamma$-invariant, dense subgroup of $\Bbb R^n$. 
Let $v_1,..,v_m\in A$ such that  $\theta:\Bbb Z^m\rightarrow A$ given by $\theta((x_i))=\sum_{i=1}^mx_iv_i$ is an isomorphism. 
 Identify the dual $\hat{A}$ of $A$ with $\hat{\Bbb Z^m}=\Bbb T^m=\Bbb R^m/\Bbb Z^m$ via the map $\hat{\theta}(\eta)=\eta\circ\theta.$

$\bullet$ Let $p:\Bbb R^n\rightarrow\Bbb R^m$ be defined by $p(a)=(\langle a,v_1\rangle,..,\langle a,v_m\rangle)$,
 where $\langle .,.\rangle$ is the usual scalar product on $\Bbb R^n$ and denote by $\pi:\Bbb R^m\rightarrow\Bbb R^m/p(\Bbb R^n)$ the
 projection.

 $\bullet$ Let $i:\Bbb T^m\rightarrow [-\frac{1}{2},\frac{1}{2})^m\subset\Bbb R^m$ be
 defined by 
$i(x+\Bbb Z^m)= (x+\Bbb Z^m)\cap [-\frac{1}{2},\frac{1}{2})^m$, for all $x\in\Bbb R^m$, and set $q=\pi\circ i:\Bbb T^m\rightarrow \Bbb R^m/p(\Bbb R^n)$.

Define $(Y,\nu)=(q(\Bbb T^m),q_*\lambda^m)$, where $\lambda^m$ is the Haar measure on $\Bbb T^m$ and $q_*\lambda^m$ is the push-forward of $\lambda^m$ through $q$. Then we have the following:

\vskip 0.05in
(1) $L^{\infty}(Y,\nu)$ is a diffuse von Neumann subalgebra of $L^{\infty}(\Bbb T^m,\lambda^m)$ (here, we consider the embedding
 $L^{\infty}(Y,\nu)\ni f\rightarrow f\circ q\in L^{\infty}(\Bbb T^m,\lambda^m)$).
\vskip 0.02in

(2) If the inclusion $A\rtimes\Gamma\rightarrow H\rtimes\Gamma$ is a resolution,
then the inclusion of von Neumann algebras $L^{\infty}(Y,\nu)\subset M:=L^{\infty}(\Bbb T^m,\lambda^m)\rtimes\Gamma$ is rigid. 
\vskip 0.02in

(3) If $p(\Bbb R^n)\cap \Bbb Z^m=\{0\}$, then $L^{\infty}(Y,\nu)'\cap M=L^{\infty}(\Bbb T^m,\lambda^m)$.
 
\endproclaim 

In the statement of this theorem we have used the fact that if  $\Gamma$ is a countable group which acts  by automorphisms on a countable abelian group $A$, then the action of $\Gamma$ on $\hat{A}$ preserves the Haar measure $h$. Also, we note that the associated crossed product von Neumann algebra $L^{\infty}(\hat{A},h)\rtimes\Gamma$ is naturally isomorphic to
the group von Neumann algebra $L(A\rtimes\Gamma)$ and that this isomorphism identifies $L^{\infty}(\hat{A},h)$ with $L(A)$.

\proclaim {2.2 Corollary} Let $n\geq 3$, $\alpha\in\Bbb N\setminus\{\beta^3|\beta\in\Bbb N\}$ and denote $\Gamma_{\alpha}=\Bbb Z[\root 3\of{\alpha}]^n\rtimes SO_n(\Bbb Z[\root 3\of{\alpha}])$. Then there exists a diffuse von Neumann subalgebra $B$ of $L(\Bbb Z[\root 3\of{\alpha}]^n)$ such that the inclusion $(B\subset L(\Gamma_{\alpha}))$ is rigid and $B'\cap L(\Gamma_{\alpha})=L(\Bbb Z[\root 3\of{\alpha}]^n)$. Moreover, if $n=3$, then $L(\Gamma_{\alpha})$ is a $\text{HT}$ factor in the sense of [Po06a, Definition 6.1.].
\endproclaim
{\it Proof.} For every $(j,k)\in S=\{1,..,n\}\times\{0,1,2\}$, let $v_{j,k}\in\Bbb R^{n}$ be given $v_{j,k}=(\alpha^{\frac{k}{3}}\delta_{i,j})_{1\leq i\leq n}.$ Then $\Bbb Z[\root 3\of{\alpha}]^n=\oplus_{(j,k)\in S}\Bbb Z v_{j,k}$.  Denote by $\langle.,.\rangle$ the natural scalar product on $\Bbb R^n$ and let $p:\Bbb R^n\rightarrow \Bbb R^{3n}=\oplus_{(j,k)\in S}\Bbb R$ be the homomorphism defined by $p(a)=(\langle a,v_{j,k}\rangle)_{j,k}$, for all $a\in\Bbb R^n$. Explicitly, if $a=(a_i)_{1\leq i\leq n}$, then $p(a)=(a_j{\alpha}^{\frac{k}{3}})_{j,k}$. Since $\alpha^{\frac{1}{3}}$ is irrational, it follows that $p(\Bbb R^n)\cap \Bbb Z^{3n}=\{0\}$.
Also, by Corollary 1.2 we have that the inclusion $\Gamma_{\alpha}\rightarrow \Bbb R^n\rtimes SO_n(\Bbb Z[\root 3\of{\alpha}])$ is a resolution. Altogether, Theorem 2.1 gives that there exists a subalgebra $B$ satisfying the conclusion.

Note that $L(\Bbb Z[\root 3\of{\alpha}]^n)$ is a Cartan subalgebra of $L(\Gamma_{\alpha})$.
Thus, in view of the first part, in order to show that $L(\Gamma_{\alpha})$ is an HT factor, it suffices to argue that SO$_3(\Bbb Z[\root 3\of{\alpha}])$ has Haagerup's property and that $\Gamma_{\alpha}$ is ICC. The first assertion is a consequence of the following general result: every countable subgroup of SO$_3(\Bbb R)$ has Haagerup's property ([GHW05], see [Co06a, Theorem 1.14]).  

To prove that $\Gamma_{\alpha}$ is ICC it suffices to show that
(1) $\{\gamma(a)-a|a\in Z[\root 3\of{\alpha}]^n\}$ is infinite, for
every $\gamma\in SO_n(\Bbb Z[\root 3\of{\alpha}])\setminus\{I\}$ and
(2) $\{\gamma(a)|\gamma\in SO_n(\Bbb Z[\root 3\of{\alpha}])\}$ is
infinite, for all $a\in Z[\root 3\of{\alpha}]^n\setminus\{0\}$. The
first assertion is clear since $\Bbb Z[\root 3\of{\alpha}]^3$ is
dense in $\Bbb R^3$. Now, since
$SO_n(\Bbb Z[\root 3\of{\alpha}])$ is an irreducible lattice in
the semisimple Lie group $SO_n(\Bbb R)\times SO_n(\Bbb C)$
and  $SO_n(\Bbb C)$ is not compact, we deduce that $SO_n(\Bbb
Z[\root 3\of{\alpha}])$ is dense in $SO_n(\Bbb R)$ (see e.g.
[Ma91] or [Wi08]). This fact  implies the second assertion.
\hfill$\blacksquare$ 

\vskip 0.1in

For the proof of Theorem 2.1 we need two technical results. To motivate and state the first result, let us  fix some notation.
 For a standard Borel space $X$ we denote by $\Cal M(X)$ the space of regular Borel probability measures on $X$ and by $B(X)$ the algebra of bounded Borel complex-valued functions on $X$. Given two measures $\mu,\nu\in\Cal M(X)$, the norm $||\mu-\nu||$ is equal to $\sup_{f\in B(X),||f||_{\infty}\leq 1}|\int_{X}fd\mu-\int_{X}fd\nu|$. 


Now, if an inclusion of the form $(A\subset A\rtimes\Gamma)$ (where $\Gamma$ is a countable group acting by automorphisms on a countable abelian group $A$) has relative property (T) then any sequence of measures $\mu_n\in\Cal M(\hat{A})$ which converge weakly to $\delta_{1}$ and are almost $\Gamma$-invariant, must ``concentrate" at the identity element $1\in\hat{A}$, i.e. $\lim_{n\rightarrow\infty}\mu_n(\{1\})=0$ ([Io09, Theorem 5.1], the converse is also true, see [Bu91]). The next proposition roughly asserts that the presence of a resolution also guarantees that almost invariant measures on $\hat{A}$  concentrate on  certain subsets.

\proclaim {2.3 Proposition} Let $H$ be a locally compact abelian group together with a dense countable subgroup $A$ and denote by $p:\hat{H}\rightarrow\hat{A}$ the map induced by restricting characters. Let $\Gamma$ be a countable group which acts by automorphisms on $H$ and leaves $A$ invariant. Suppose that the inclusion $A\rtimes\Gamma \rightarrow H\rtimes\Gamma$ is a resolution. Also, let $V\subset \hat{H}$ be a $\Gamma$-invariant  neighborhood of $1\in \hat{H}$.

Then for any sequence of measures $\mu_n\in\Cal M(\hat{A})$ which converge weakly to $\delta_{1}$ and satisfy $\lim_{n\rightarrow\infty}||\gamma_{*}\mu_n-\mu_n||=0$, for all $\gamma\in\Gamma$, we have that $\lim_{n\rightarrow\infty}\mu_n(p(V))=1$. In particular, we have that $\lim_{n\rightarrow\infty}\mu_n(p(\hat{H}))=1$.
\endproclaim
{\it Proof.} 
Given $V\subset\hat{H}$ and a sequence $\{\mu_n\}_{n\geq 1}\subset\Cal M(\hat{A})$ as in the hypothesis, we begin by showing:
\vskip 0.05in
\noindent
{\it Claim. } There exists $n$ such that $\mu_n(p(V))>0$. 
\vskip 0.05in
\noindent
{\it Proof of the claim.} 
Let us first prove the claim under the additional assumption that $\mu_n$ is $\Gamma$-quasi-invariant, for all $n$.
Fix $n\geq 1$. Since $\mu_n$ is $\Gamma$-quasi-invariant, we can define $g_{\gamma}=(d(\gamma_*\mu_n)/d\mu_n)^{\frac{1}{2}}$, for all $\gamma\in\Gamma$, where $d(\gamma_*\mu_n)/d\mu_n$ denotes the Radon-Nikodym derivative of $\gamma_*\mu_n$ with respect to $\mu_n$. Next, we see every $a\in A$ as a character on $\hat{A}$ and therefore as as a function in $L^{\infty}(\hat{A},\mu_n)$. Then the  formulas $$\pi_n(a)(f)=af,\pi_n(\gamma)(f)=g_{\gamma}(f\circ{\gamma}^{-1}),$$ for all $a\in A,\gamma\in\Gamma$ and $f\in L^2({\hat{A},\mu_n})$,
define a unitary representation $\pi_n:A\rtimes\Gamma\rightarrow\Cal U(L^2(\hat{A},\mu_n))$. Let $\xi_n=1_{\hat{A}}\in L^2(\hat{A},\nu_n)$. For all $a\in A$ and $\gamma\in\Gamma$, we have that $$||\pi_n(\gamma)(\xi_n)-\xi_n||=||g_{\gamma}-1||_2\leq ||g_{\gamma}^2-1||_1^{\frac{1}{2}}=||\gamma_*\mu_n-\mu_n||^{\frac{1}{2}},$$
$$||\pi_n(a)(\xi_n)-\xi_n||=(\int_{\hat{A}}|\eta(a)-1|^2d\mu_n(\eta))^{\frac{1}{2}}.$$

Using the assumptions made on $\mu_n$, it follows that the vectors $\xi_n$ form an almost invariant sequence for the representation $\pi=\oplus_{n\geq 1}\pi_n:A\rtimes\Gamma\rightarrow \Cal U(\Cal H)$, where $\Cal H=\oplus_{n\geq 1}L^2(\hat{A},\mu_n)$. 
Since the inclusion $A\rtimes\Gamma\rightarrow H\rtimes\Gamma$ is a resolution, we can find a $\pi$-invariant Hilbert subspace $\Cal K\subset\Cal H$ and a unitary representation $\sigma:H\rtimes\Gamma\rightarrow\Cal U(\Cal K)$ which has almost invariant vectors and satisfies $\sigma(g)=\pi(g)_{|\Cal K}$, for all $g\in A\rtimes\Gamma$. 
Let $\{\zeta_k\}_{k\geq 1}\subset \Cal K$ be a sequence of $\sigma$-almost invariant unit vectors. For every $k$, let $\nu_k\in \Cal M(\hat{H})$ be given by $\langle\sigma(g) \zeta_k,\zeta_k\rangle=\int_{\hat{H}}\eta(g)d\nu_k(\eta),$ for all $g\in H$. Notice that $\nu_k$ converge weakly to $\delta_1$, as $k\rightarrow\infty$.

Next, if we  set $\rho_k=p_*\nu_k\in\Cal M(\hat{A})$, then for  each $a\in A$ we have that $$\int_{\hat{A}}a d\rho_k=\int_{\hat{H}}(a\circ p) d\nu_k=\int_{\hat{H}}\eta(a)d\nu_k(\eta)=\tag a$$  $$\langle\sigma(a)\zeta_k,\zeta_k\rangle=\langle\pi(a)(\zeta_k),\zeta_k\rangle. $$
Now, for every $k\geq 1$, decompose $\zeta_k=\sum_{n\geq 1}\zeta_k^n$, where $\zeta_k^n\in L^2(\hat{A},\mu_n)$. Thus, for all $a\in A$, we have that  $$\langle\pi(a)(\zeta_k),\zeta_k\rangle=\sum_{n\geq 1}\int_{\hat{A}}a|\zeta_k^n|^2d\mu_n\tag b$$ By combining (a) and (b) we deduce that $d\rho_k=\sum_{n\geq 1}|\zeta_k^n|^2d\mu_n$, for all $k\geq 1$. Since $\rho_k(p(V))=\nu_k(V)$ and $\nu_k\rightarrow\delta_{1}$ weakly, we get that $\lim_{k\rightarrow\infty}\rho_k(p(V))=1$. Thus, we can find $n$ such that $\mu_n(p(V))>0$.

In general, if $\mu_n$ are not necessarily quasi-invariant,  let $\{\gamma_i\}_{i\geq 1}$ be an enumeration of $\Gamma$.   For every $n$, set $\mu_n'=\sum_{i\geq 1}\frac{1}{2^i}{\gamma_i}_*\mu_n$. Then $\mu_n'$ are $\Gamma$ quasi-invariant measures which satisfy the hypothesis. By applying the first part of the proof, we get that $\mu_n(p(V))=\mu_n'(p(V))>0$, for some $n$. 
\hfill$
\square$
\vskip 0.1in
Suppose by contradiction that the conclusion of the theorem is false. Then, after passing to a subsequence, we can assume that  $\lim_{n\rightarrow\infty}\mu_n(p(V))=c<1.$ Thus, for large enough $n$ we have that $\mu_n(\hat{A}\setminus p(V))>0$, so can define $\mu_n'\in\Cal M(\hat{A})$ by letting $\mu_n'(X)=\frac{\mu_n(X\setminus p(V))}{\mu_n(\hat{A}\setminus p(V))}$, for every Borel set $X\subset \hat{A}$.  Notice that $\mu_n'\rightarrow\delta_1$ weakly, as $n\rightarrow\infty$.  To see this, just remark that for every neighborhood $W$ of $1\in\hat{A}$ we have that $\mu_n'(\hat{A}\setminus W)\leq \frac{\mu_n(\hat{A}\setminus W)}{\mu_n(\hat{A}\setminus p(V))}\rightarrow 0$, as $n\rightarrow\infty$.

Next, it is easy to see that since $V$ is $\Gamma$-invariant, we get that $\lim_{n\rightarrow\infty}||\gamma_*\mu_n'-\mu_n'||=0$, for all $\gamma\in\Gamma$. Altogether, it follows that $\mu_n'$ satisfy the conditions of the hypothesis. Therefore, we can apply the claim  and derive that $\mu_n'(p(V))>0$, for some $n$, a contradiction.\hfill$\blacksquare$

\vskip 0.05in 

The second ingredient needed in the proof of Theorem 2.1 is the following criterion for rigidity which we derive as a consequence of results from [Io09].

\proclaim {2.4 Proposition} Let $\Gamma\curvearrowright (X,\mu)$ be a measure preserving action of a countable group $\Gamma$ on a standard
probability space $(X,\mu)$. Let $p^i:X\times X\rightarrow X$ be the projection $p^i(x_1,x_2)=x_i$, for $i\in\{1,2\}$, and endow $X\times X$ 
with the diagonal action of $\Gamma$. 
Let $(Y,\nu)$ be another probability space together with a measurable, measure preserving onto 
map $q:X\rightarrow Y$. View $L^{\infty}(Y,\nu)$ as a von Neumann subalgebra of $L^{\infty}(X,\mu)$ via the 
embedding $L^{\infty}(Y,\nu)\ni f\rightarrow f\circ q\in L^{\infty}(X,\mu)$.
\vskip 0.03in
Assume that for any sequence of measures $\{\nu_n\}_{n\geq 1}\subset\Cal M(X\times X)$ such that $p^i_*(\nu_n)=\mu$, for all $i\in\{1,2\}$
 and $n\geq 1$, 

$(i)$ $\lim_{n\rightarrow\infty}\int_{X\times X}f_1(x)f_2(y)d\nu_n(x,y)=\int_{X}f_1f_2d\mu,$ for all $f_1,f_2\in B(X)$, and

$(ii)$ $\lim_{n\rightarrow\infty}||\gamma_*\nu_n-\nu_n||=0$, for all $\gamma\in\Gamma$,

we have that $\lim_{n\rightarrow\infty}\nu_n(\{(x,y)\in X\times X|q(x)=q(y)\})=1$.

 Then the inclusion of von Neumann algebras 
$L^{\infty}(Y,\nu)\subset  L^{\infty}(X,\mu)\rtimes\Gamma$ is rigid.
\endproclaim

{\it Proof.} Denote $M=L^{\infty}(X,\mu)\rtimes\Gamma$ and let $(\Cal H_n,\xi_n)$ be a sequence of pointed Hilbert $M$-bimodules such that $\lim_{n\rightarrow\infty}||z\xi_n-\xi_n z||=0$, for all $z\in M$.  To get the conclusion we have to show that there exists a sequence $\eta_n\in\Cal H_n$ of $L^{\infty}(Y,\nu)$-central vectors such that $\lim_{n\rightarrow\infty}||\eta_n-\xi_n||=0$ (see Definition 1.4). 

By  [Io09, Lemma 2.1] we can find a sequence $\{\nu_n\}_{n\geq 1}\subset\Cal M(X\times X)$
 which verifies all the conditions from the hypothesis and  satisfies $\int_{X\times X}f_1(x)f_2(y)d\nu_n(x,y)=\langle f_1\xi_nf_2,\xi_n\rangle,$ for all $f_1,f_2\in B(X)$. Thus, if $\Delta_q:=\{(x,y)\in X\times X|q(x)=q(y)\}$, then $\lim_{n\rightarrow\infty}\nu_n(\Delta_q)=1$.

Next, for every $f_1,f_2\in B(X)$, let $f_1\otimes f_2:X\rightarrow\Bbb C$ be given by $(f_1\otimes f_2)(x,y)=f_1(x)f_2(y)$. Notice then that by the way $\nu_n$ is defined, the map $L^2(X\times X,\nu_n)\ni f_1\otimes f_2\rightarrow f_1\xi_n f_2\in\Cal H_n$, for every $f_1,f_2\in B(X)$, extends to an embedding of Hilbert $L^{\infty}(X,\mu)$-bimodules $\theta_n:L^2(X\times X,\nu_n)\rightarrow\Cal H_n$.  Here, on $L^2(X\times X,\nu_n),$ we consider the $L^{\infty}(X,\mu)$-bimodule structure given by $f_1\cdot g\cdot f_2=(f_1\otimes f_2)g$, for all $f_1,f_2\in L^{\infty}(X,\mu)$ and $g\in L^2(X\times X,\nu_n)$.
Let $\eta_n=\theta_n(1_{\Delta_q})$. Since $1_{\Delta_q}\in L^2(X\times X,\nu_n)$ is an $L^{\infty}(Y,\nu)$-central vector, we get that $\eta_n\in\Cal H_n$ is an $L^{\infty}(Y,\nu)$-central vector. Finally, notice that $||\eta_n-\xi_n||=||1_{\Delta_{q}}-1_{X\times X}||_{L^2(X\times X,\nu_n)}=\sqrt{\nu_n((X\times X)\setminus\Delta_q)}\rightarrow 0$, as $n\rightarrow\infty$.\hfill$\blacksquare$

\vskip 0.1in
We are now ready to prove  2.1:
\vskip 0.05in
\noindent
{\it Proof of Theorem 2.1}. (1) To derive that $L^{\infty}(Y,\nu)$ is  diffuse, we only need to show that $\lambda^m(q^{-1}(\{y\}))=0$, for every $y\in Y$.
This is clear since $n<m$ and $q^{-1}(\{y\})\subset (y+p(\Bbb R^n))+\Bbb Z^m\subset\Bbb T^m$. 
\vskip 0.1in

\noindent
(2) To prove the rigidity assertion, 
let $\nu_k\in\Cal M(\Bbb T^m\times\Bbb T^m)$ be a sequence of measures such that 
 $$\lim_{k\rightarrow\infty}\int_{\Bbb T^m\times \Bbb T^m}f_1(x)f_2(y)d\nu_k(x,y)=\int_{\Bbb T^m}f_1f_2d\lambda^m,\forall f_1,f_2\in B(\Bbb T^m)\tag a$$ and
 $$\lim_{k\rightarrow\infty}||\gamma_*\nu_k-\nu_k||=0,\forall\gamma\in\Gamma\tag b$$
Denote by $\Delta_q=\{(x,y)\in\Bbb T^m\times \Bbb T^m|q(x)=q(y)\}$. Following Proposition 2.4 in order to get the conclusion it 
suffices to argue that $\lim_{k\rightarrow\infty}\nu_k(\Delta_q)=1.$ To this end, notice that (a) gives that for every bounded Borel function $f$ on $\Bbb T^m$ we have that $\lim_{k\rightarrow\infty}\int_{\Bbb T^m\times\Bbb T^m}|f(x)-f(y)|^2d\nu_k(x,y)=0$. This implies that $\lim_{k\rightarrow\infty}\int_{\Bbb T^m\times\Bbb T^m}||i(x)-i(y)||^2d\nu_k(x,y)=0,$ where $||.||$ denotes the Euclidian norm on $\Bbb R^m$. Thus, we deduce that $$\lim_{k\rightarrow\infty}\nu_k(\{(x,y)\in\Bbb T^m\times\Bbb T^m|||i(x)-i(y)||\leq \frac{1}{2}\})=1\tag c$$

Next, assume that $x,y\in\Bbb T^m$ satisfy $||i(x)-i(y)||\leq \frac{1}{2}$. Since $i(x)-i(y)\in (x-y)+\Bbb Z^m$, we deduce that $i(x)-i(y)=i(x-y)$ and therefore that $q(x)-q(y)=q(x-y)$. By combining this fact with (c), we get that $\lim_{k\rightarrow\infty}\nu_k(\{(x,y)\in\Bbb T^m\times\Bbb T^m|q(x)-q(y)=q(x-y)\})=1$. Thus, showing that $\lim_{k\rightarrow\infty}\nu_k(\Delta_q)=1$ is equivalent to proving that $$\lim_{k\rightarrow\infty}\nu_k(\{(x,y)\in\Bbb T^m\times\Bbb T^m|q(x-y)=0\})=1\tag d$$

If we let $r:\Bbb T^m\times\Bbb T^m\rightarrow\Bbb T^m$ be given by $r(x,y)=x-y$ and for every $k$, define $\mu_k=r_*\nu_k\in\Cal M(\Bbb T^m)$,
then (d) can be rewritten as $\lim_{k\rightarrow\infty}\mu_k(\{x\in\Bbb T^m|q(x)=0\})=1$.

Now, note that the inclusion $A=\Bbb Zv_1\oplus...\oplus\Bbb Zv_m\subset H=\Bbb R^n$ gives rise to a homomorphism $p':\hat{H}\simeq\Bbb R^n\rightarrow\hat{A}\simeq\Bbb T^m$ given by $p'(a)=(\langle a,v_1\rangle+\Bbb Z,..,\langle a,v_m\rangle+\Bbb Z)$, for all $a\in\Bbb R^n$. In other words,  $p'$ is the composition between $p$ and the projection $\Bbb R^m\rightarrow\Bbb T^m$. 
Let $\varepsilon>0$ such that $V=\{a\in\Bbb R^n|||a||\leq \varepsilon\}$ satisfies $p(V)\subset [-\frac{1}{2},\frac{1}{2})^{m}$.
We claim that $p'(V)\subset \{x\in\Bbb T^m|q(x)=0\}$. Indeed, if $a\in V$, then  $p(a)\in  [-\frac{1}{2},\frac{1}{2})^{m}$, thus $i(p'(a))=p(a)$, hence $q(p'(a))=(\pi\circ i)(p'(a))=\pi(p(a))=0$.

On the other hand, $V$ is a $\Gamma$-invariant neighborhood of $0\in\Bbb R^n$. 
Also, remark that (a) implies that $\mu_k$ converge weakly to $\delta_{0}$ while (b) implies that $\lim_{k\rightarrow\infty}||\gamma_*\mu_k-\mu_k||=0$, for all $\gamma\in\Gamma$.  By applying Proposition 2.3 we deduce that $\lim_{k\rightarrow\infty}\mu_k(p'(V))=1$. Since  $p'(V)\subset \{x\in\Bbb T^m|q(x)=0\}$, we get that $\lim_{k\rightarrow\infty}\mu_k(\{x\in\Bbb T^m|q(x)=0\})=1$. This proves (d) and thus the conclusion.
\vskip 0.1in

\noindent
(3) Firstly, it is easy to see that the conclusion is equivalent to $\lambda^m(\{x\in\Bbb T^m|q(\gamma x)=q(x)\})=0$, for all $\gamma\in\Gamma\setminus\{I\}$ (where $I$ denotes the identity matrix). Assuming that this is not the case we can find $\gamma\in\Gamma\setminus\{I\}$ such that $\lambda^m(\{x\in\Bbb T^m|q(\gamma x)=q(x)\})>0$. By using the definition of $q$, this implies that $$\lambda^m(\{x\in\Bbb T^m|\gamma x-x\in p(\Bbb R^n)+\Bbb Z^m\})>0\tag e$$ Secondly, notice that the action of $\Gamma$ on $\Bbb Z^m\simeq A$ is realized through a homomorphism $\rho:\Gamma\rightarrow$ GL$_m(\Bbb Z)$. The (dual) action of $\Gamma$ on $\Bbb T^m\simeq\hat{A}$ is then given by $\gamma (x+\Bbb Z^m)=(\rho(\gamma)^{-1})^t(x)+\Bbb Z^m$. 
Altogether, (e) implies that if $\mu^m$ denotes the Lebesgue measure on $\Bbb R^m$, then $\mu^m(\{x\in\Bbb R^m| (\rho(\gamma)^{-1})^t(x)-(x)\in p(\Bbb R^n)+\Bbb Z^m\})>0.$ 
This easily implies that $$(\rho(\gamma^{-1})^t-I)(\Bbb R^m)\subset p(\Bbb R^n)\tag f$$

Finally, since $\rho(\gamma^{-1})\in$ SL$_m(\Bbb Z)$, we get that $(\rho(\gamma^{-1})^t-I)(\Bbb Z^m)\subset \Bbb Z^m$. By combining this with (f) and the fact that $p(\Bbb R^n)\cap\Bbb Z^m=\{0\}$, we deduce that $(\rho(\gamma^{-1})^t-I)(\Bbb Z^m)=\{0\}$. This means that $\rho(\gamma)=I$, or, equivalently, that $\gamma$ acts trivially on $A$. Since $A$ is dense in $\Bbb R^n$, we further get that $\gamma$ acts trivially on $\Bbb R^n$ which implies that $\gamma=I$, a contradiction.
\hfill$\blacksquare$
\vskip 0.2in
\head 3. Relative property (T) subsets of semidirect product groups.\endhead
\vskip 0.1in
In this section we show that Haagerup's property is not preserved under generalized wreath products.  Using this fact we give  the first examples of von Neumann algebras which neither have Haagerup's property nor admit any diffuse rigid  von Neumann subalgebras. We start with the following result which asserts that if a group $\Gamma$ acts on an abelian group $A$ through a quotient group $\Gamma_0$, then the presence of a relative property (T) subset in $\Gamma_0$ (or the lack of Haagerup's property)   is inherited by the semidirect product $A\rtimes\Gamma$.

\proclaim {3.1 Theorem} Let $\Gamma_0$ be a countable group and let $\rho:\Gamma_0\rightarrow\text{Aut}(A)$ be an action by automorphisms on a countable abelian group $A$. Let $\Gamma$ be a countable group together with a surjective homomorphism $p:\Gamma\rightarrow\Gamma_0$ and consider the action of $\Gamma$ on $A$ given by $\tilde\rho=\rho\circ p:\Gamma\rightarrow$Aut$(A)$. 
\vskip 0.03in
(1) Suppose that $X$ is a subset of $\Gamma_0$ such that the inclusion $(X\subset\Gamma_0)$ has relative property (T). For $a\in A$, let $X_a=\{\rho(\gamma)(a)|\gamma\in X\}\subset A$. Then the semidirect product $A\rtimes_{\tilde\rho}\Gamma$ has relative property (T) with respect to $X_a$, for every $a\in A$.
\vskip 0.03in
(2) Assume that there is $a\in A$ such that its stabilizer $\{\gamma\in\Gamma_0|\rho(\gamma)(a)=a\}$ in $\Gamma_0$ is finite.
If $\Gamma_0$ does not have Haagerup's property then $A\rtimes_{\tilde\rho}\Gamma$  does not have Haagerup's property.
\endproclaim
{\it Proof.} (1) Fix $a\in A$ and let $\pi:A\rtimes_{\tilde\rho}\Gamma\rightarrow\Cal U(\Cal H)$ be a unitary representation which admits a sequence $\{\xi_n\}_{n\geq 1}\subset\Cal H$ of almost invariant, unit vectors. To get the conclusion we have to show that $\xi_n$ are uniformly $\pi(X_a)$-almost invariant, i.e. $\lim_{n\rightarrow\infty}\sup_{\gamma\in X}||\pi(\rho(\gamma)(a))(\xi_n)-\xi_n||=0$.

 Firstly, for every $n\geq 1$, let $\mu_n\in\Cal M(\hat{A})$ such that $\langle\pi(a)\xi_n,\xi_n\rangle=\int_{\hat{A}}ad\mu_n$, for each $a\in A$. By the proof of [Bu91, Proposition 7] we have that $||\gamma_*\mu_n-\mu_n||\leq 2||\pi(\gamma)(\xi_n)-\xi_n||$,  for all $\gamma\in\Gamma.$ Here, on $\hat{A}$ we consider the natural actions of $\Gamma$ and $\Gamma_0$ induced by $\tilde\rho$ and $\rho$, respectively. For every $\gamma\in\Gamma_0$, fix $\tilde\gamma\in\Gamma$ such that $p(\tilde\gamma)=\gamma$. Then the above implies that $$||\gamma_*\mu_n-\mu_n||=||{\tilde\gamma}_*\mu_n-\mu_n||\leq 2||\pi(\tilde\gamma)(\xi_n)-\xi_n||,\forall\gamma\in\Gamma_0\tag a$$

Secondly, let $\{\gamma_i\}_{i\geq 1}$ be an enumeration of $\Gamma_0\setminus\{e\}$. For every $n\geq 1$, let $\nu_n\in\Cal M(\hat{A})$ be given by $\nu_n=(1-\frac{1}{2^n})\mu_n+\sum_{i\geq 1}\frac{1}{2^{i+n}}{\gamma_i}_*\mu_n$. Then we have that $||\nu_n-\mu_n||\leq \frac{1}{2^{n-1}}$, for all $n\geq 1$, and thus (a) implies that $$||\gamma_*\nu_n-\nu_n||\leq 2||\nu_n-\mu_n||+||\gamma_*\mu_n-\mu_n||\leq\tag b$$ $$ \frac{1}{2^{n-2}}+2||\pi(\tilde\gamma)(\xi_n)-\xi_n||,\forall\gamma\in\Gamma_0.$$

Next, fix $n\geq 1$. Notice that $\nu_n$ is a $\Gamma_0$-quasi-invariant measure and let $g_{\gamma}=(d(\gamma_*\nu_n)/d\nu_n)^{\frac{1}{2}}\in L^2(\hat{A},\nu_n)$, for every $\gamma\in\Gamma$. The formula $\sigma_n(\gamma)(f)=g_{\gamma}(f\circ{\gamma}^{-1}),$ for all $f\in L^2(\hat{A},\nu_n)$ and $\gamma\in\Gamma$, defines a unitary representation $\sigma_n:\Gamma_0\rightarrow\Cal U(L^2(\hat{A},\nu_n))$. If $\eta_n=1_{\hat{A}}\in L^2(\hat{A},\nu_n)$, then, as in the proof of Proposition 2.3, we have that $$||\sigma_n(\gamma)(\eta_n)-\eta_n||\leq ||\gamma_*\nu_n-\nu_n||^{\frac{1}{2}},\forall\gamma\in\Gamma_0\tag c$$ 

Since the vectors $\xi_n$ are $\pi(\Gamma)$-almost invariant, by combining (b) and (c) we deduce that $\lim_{n\rightarrow\infty}||\sigma_n(\gamma)(\eta_n)-\eta_n||=0$, for each $\gamma\in\Gamma_0$. 
Since the inclusion $(X\subset\Gamma_0)$ has relative property (T), by [Co06b, Theorem 1.1] we get that $\varepsilon_n:=\sup_{\gamma\in X}||\sigma_n(\gamma)(\eta_n)-\eta_n||\rightarrow 0,$ as $n\rightarrow\infty$. Now, for $\gamma\in X$ we have that $$||\gamma_{*}\nu_n-\nu_n||=||g_{\gamma}^2-1||_1\leq ||g_{\gamma}+1||_2||g_{\gamma}-1||_2\leq 2||g_{\gamma}-1||_2=$$ $$2||\sigma_n(\gamma)(\eta_n)-\eta_n||\leq 2\varepsilon_n,\forall n\geq 1.$$ Thus, we get that $||\gamma_*\mu_n-\mu_n||\leq 2||\nu_n-\mu_n||+||\gamma_*\nu_n-\nu_n||\leq \frac{1}{2^{n-1}}+2\varepsilon_n,$ for all $\gamma\in X.$ This implies that for every $\gamma\in X$ we have that $$|\langle\pi(\rho(\gamma)(a))(\xi_n),\xi_n\rangle-\langle\pi(a)(\xi_n),\xi_n\rangle|=|\int_{\hat{A}}\rho(\gamma)(a) d\mu_n-\int_{\hat{A}}a d\mu_n|=\tag d$$ $$|\int_{\hat{A}}(a\circ{\gamma^{-1}}) d\mu_n-\int_{\hat{A}}a d\mu_n|\leq ||\gamma_*\mu_n-\mu_n||\leq \frac{1}{2^{n-1}}+2\varepsilon_n.$$

Finally, (d) together with a standard calculation gives that $$||\pi(\rho(\gamma)(a))(\xi_n)-\xi_n||^2\leq ||\pi(a)(\xi_n)-\xi_n||^2+2( \frac{1}{2^{n-1}}+2\varepsilon_n),\forall\gamma\in X,$$ which proves the conclusion.
\vskip 0.05in
\noindent
(2) Assume by contradiction that $A\rtimes_{\tilde\rho}\Gamma$ has Haagerup's property while $\Gamma_0$ does not have it. Thus we can find a $c_0$ unitary representation $\pi:A\rtimes_{\tilde\rho}\Gamma\rightarrow\Cal U(\Cal H)$ which admits a sequence $\{\xi_n\}_{n\geq 1}\subset\Cal H$ of almost invariant, unit vectors.  Let $\sigma_n:\Gamma_0\rightarrow\Cal U(L^2(\hat{A},\nu_n))$ and $\eta_n\in L^2(\hat{A},\nu_n)$ be constructed as in the proof of part (1). Recall that $\eta_n$ are almost invariant unit vectors, i.e. $\lim_{n\rightarrow\infty}||\sigma_n(\gamma)(\eta_n)-\eta_n||=0$, for each $\gamma\in\Gamma_0$. 

Since $\Gamma_0$ does not have Haagerup's property  by [Pe09, Theorem 2.6.] we can find an infinite subset $X$ of $\Gamma_0$ and an increasing sequence $\{k_n\}_{n\geq 1}$ of natural numbers such that $\lim_{n\rightarrow\infty}\sup_{\gamma\in X}||\sigma_{k_n}(\gamma)(\eta_{k_n})-\eta_{k_n}||=0$. Let $a\in A$ such that its stabilizer in $\Gamma_0$ is finite. The last part of the proof of (1) implies that $\lim_{n\rightarrow\infty}\sup_{\gamma\in X}||\pi(\rho(\gamma)(a))(\xi_{k_n})-\xi_{k_n}||=0$. On the other hand, since the stabilizer of $a$ in $\Gamma_0$ is finite and $\pi$ is a $c_0$ representation, we get that $\lim_{\gamma\rightarrow\infty}\langle\pi(\rho(\gamma)(a))(\xi_n),\xi_n\rangle=0$. Altogether, this gives a contradiction as $X$ is infinite.
 \hfill$\blacksquare$
\vskip 0.1in

\noindent
{\bf 3.2 Remarks.} (1) We note that the proof of part (1) in fact shows more: if $X\subset A$ is a set such that the inclusion $(X\subset A\rtimes_{\rho}\Gamma_0)$ has relative property (T), then the inclusion $(X\subset A\rtimes_{\tilde\rho}\Gamma)$ has relative property (T).
\vskip 0.03in
\noindent
(2) Let us also remark that the proof of (2) can be adapted to show that if $A\rtimes_{\rho}\Gamma_0$ is not Haagerup then $A\rtimes_{\tilde\rho}\Gamma$ is not Haagerup, provided that the stabilizer of some $a\in A$ in $\Gamma_0$ is finite.  Indeed, in the notations from above, if $A\rtimes_{\rho}\Gamma_0$ is not Haagerup, then we can find an infinite set $X\subset A\rtimes_{\rho}\Gamma_0$ and a sequence $\{k_n\}_{n\geq 1}$ such that  $\lim_{n\rightarrow\infty}\sup_{\gamma\in X}||\sigma_{k_n}(\gamma)(\eta_{k_n})-\eta_{k_n}||=0$. Thus, $\lim_{n\rightarrow\infty}\sup_{\gamma\in X}||\sigma_{k_n}(\gamma a\gamma^{-1})(\eta_{k_n})-\eta_{k_n}||=0$. If the projection of $X$ onto $\Gamma_0$ is infinite, then  $\{\gamma a\gamma^{-1}|\gamma\in X\}$ is an infinite subset of $A$ and a contradiction is reached as in the above proof. Otherwise,  the set $Y$ of all $a\in A$ such that $(a,\gamma)\in X$, for some $\gamma\in \Gamma_0$, is infinite. Since the projection of $X$ onto $\Gamma_0$ is finite it is clear that $\lim_{n\rightarrow\infty}\sup_{a\in Y}||\sigma_{k_n}(a)(\eta_{k_n})-\eta_{k_n}||=0$. Again, we obtain a contradiction as in the end of the proof of part (1). 
\vskip 0.1in

Recently, Y. de Cornulier, Y. Stalder and A. Valette have shown that if two countable groups $A$ and $\Gamma$  have Haagerup's property then so does their {\it wreath product} $A\wr \Gamma=A^{\Gamma}\rtimes\Gamma$ ([CSV09]). More generally, they consider {\it generalized wreath product} groups  $A\wr_{X}\Gamma=A^{X}\rtimes\Gamma$, where $X$ is a countable $\Gamma$-set and $\Gamma$ acts on $A^{X}=\bigoplus_{x\in X}A$ by shifting indices. In this context, they show that if $A$ and $\Gamma$ are Haagerup then so are certain  wreath products $A\wr_{X}\Gamma$ (e.g. when $X$ is a Haagerup  quotient group of $\Gamma$) ([CSV09, Theorem 6.2.]). Furthermore, it is conjectured in [CSV09] that all such wreath products are Haagerup. As a consequence of Theorem 3.1, we disprove this conjecture by showing, for example, that if $X$ is an infinite quotient of $\Gamma$ with property (T), then $A\wr_{X}\Gamma$ is not Haagerup. 

\proclaim {3.3 Corollary} Let $A$ be a non-trivial countable group. Let $\Gamma$ be a countable group together with a quotient group $\Gamma_0$. Endow $\Gamma_0$  with the left multiplication action of $\Gamma$. Then the generalized wreath product $A\wr_{\Gamma_0}\Gamma$ is Haagerup if and only if  $A,\Gamma$ and $\Gamma_0$ are.
\endproclaim
{\it Proof.} The {\it if} part is a particular case of [CSV09, Theorem 6.2]. To prove the {\it only if} part, assume that $A\wr_{\Gamma_0}\Gamma$ is Haagerup. Since both $A$ and $\Gamma$ are subgroups of $A\wr_{\Gamma_0}\Gamma$, they must be Haagerup. Let $a\in A\setminus\{e\}$ and let $A_0$ be the cyclic group generated by $a$. Since $A_0\wr_{\Gamma_0}\Gamma$ is Haagerup (being a subgroup of $A\wr_{\Gamma_0}\Gamma$) and $A_0$ is abelian, Theorem 3.1 (b) implies that $\Gamma_0$ is Haagerup.\hfill$\blacksquare$
\vskip 0.1in

The above corollary gives new examples of countable groups which are not Haagerup and yet do not admit any infinite subgroups with relative property (T).  More precisely, in the above context, $G=A\wr_{\Gamma_0}\Gamma$ is such a group, whenever $A$, $\Gamma$ are Haagerup and $\Gamma_0$ is not. Indeed, by Corollary 3.3, $G$ is not Haagerup, while by [CSV09, Theorem 6.7.], $G$ does not have relative property (T) with respect to any infinite subgroup. Moreover, we can use results from [Po06ab] and [Io07] to deduce that the group von Neumann algebra of $G$ has no rigid von Neumann subalgebra: 

\proclaim {3.4 Corollary} Let $A$ be a non-trivial countable Haagerup group. Let $\Gamma$ be a countable Haagerup together with a quotient group $\Gamma_0$. Assume that $\Gamma_0$ is not Haagerup and endow it with the left multiplication action of $\Gamma$. Then the group von Neumann algebra $N=L(A\wr_{\Gamma_0}\Gamma)$ does not have Haagerup property and does not admit any diffuse von Neumann subalgebra $B$ such that the inclusion $(B\subset  N)$ is rigid.
\endproclaim
{\it Proof.} By Corollary 3.3, $A\wr_{\Gamma_0}\Gamma$ does not have Haagerup's property and therefore its group von Neumann algebra does not have it either (see e.g. [Po06a]). Now, assume by contradiction that the inclusion $(B\subset  N)$ is rigid, for some diffuse von Neumann subalgebra $B$ of $N$. Since $\Gamma$ has Haagerup's property, the proof of [Po06a, Theorem 6.2.] gives that a corner of $B$ embeds into $L(A^{\Gamma_0})$, in the sense of [Po06b, Section 2].

To get a contradiction, we apply results from [Io07], leaving the details to the reader. 
Notice first that $N$ can be written as $(\overline{\bigotimes}_{\Gamma_0}L(A))\rtimes\Gamma$, where $\Gamma$ acts on $\overline{\bigotimes}_{\Gamma_0}L(A)$ by Bernoulli shifts. Using this observation and the rigidity of the inclusion $(B\subset N)$, the proof of [Io07, Theorem 3.6.] implies that a corner of $B$ can be embedded into $L(A^{F})$, for some finite subset $F$ of $\Gamma_0$. Finally, the proof of [Io07, Corollary 3.7.] shows that since $L(A)$ has Haagerup's property, $B$ cannot be diffuse, a contradiction.\hfill$\blacksquare$

\head  References\endhead
\item {[Bu91]} M. Burger: {\it Kazhdan constants for SL(3,Z),} J. Reine Angew. Math. {\bf 413} (1991),
36--67.
\item {[Ch83]} M. Choda: {\it Group factors of the Haagerup type,} Proc. Japan Acad., {\bf 59} (1983), 174--177.
\item  {[CJ85]} A. Connes, V.F.R. Jones: {\it Property (T) for von Neumann algebras}, Bull. London Math.
Soc. {\bf 17} (1985), 57--62.
\item {[CCJJV01]} P.A. Cherix, M. Cowling, P. Jolissaint, P. Julg, A. Valette: {\it Groups with the
Haagerup Property}, Birkh¨auser, Progress in Mathematics 197, 2001.
\item {[Co06a]} Y. de Cornulier: {\it Kazhdan and Haagerup Properties in algebraic groups over local fields,} J. Lie Theory {\bf 16} (2006), 67--82. 
\item {[Co06b]} Y. de Cornulier: {\it Relative Kazhdan Property}, Ann. Sci. Ecole Norm. Sup. {\bf 39}  (2), (2006), 301--333. 
\item {[CSV09]} Y. de Cornulier, Y. Stalder, A. Valette: {\it Proper actions of wreath products and generalizations,} preprint arXiv:0905.3960. 
\item {[Fu09]} A. Furman: {\it A survey of Measured Group Theory}, preprint arXiv:0901.0678.
\item {[Io07]} A. Ioana: {\it Rigidity results for wreath product II$_1$ factors}, Journal of Functional Analysis {\bf 252} (2007) 763--791.
\item {[Io09]} A. Ioana: {\it Relative Property (T) for the Subequivalence Relations Induced by the Action of SL$_2(\Bbb Z)$ on $\Bbb T^2$}, 
preprint arXiv:0901.1874.
\item {[Ka67]} D. Kazhdan: {\it On the connection of the dual space of a group with the structure of
its closed subgroups}, Funct. Anal. and its Appl. {\bf 1} (1967), 63-–65. 
\item {[Ma82]} G. Margulis: {\it Finitely-additive invariant measures on Euclidian spaces}, Ergodic Theory
Dynam. Systems {\bf 2} (1982), 383–396.
\item {[Ma91]} G. Margulis: {\it Discrete subgroups of semisimple Lie groups}, Springer, 1991.
\item {[MvN36]} F.J. Murray, J. Von Neumann: {\it On rings of operators}, Ann. of Math. (2) {\bf 37} (1936),
no. 1, 116–-229.
\item {[Pe09]} J. Peterson: {\it Examples of group actions which are virtually W*-superrigid}, preprint 2009.
\item {[Po06a]} S. Popa: {\it On a class of type II1 factors with Betti numbers invariants,} Ann. of Math. (2) {\bf 163} (2006), no. 3, 809-–899.
\item {[Po06b]} S. Popa: {\it Strong rigidity of II$_1$ factors arising from malleable actions of w-rigid groups
 I}, Invent. Math. {\bf 165} (2006), 369--408. 
\item {[Po07]} S. Popa: {\it Deformation and rigidity for group 
actions and von Neumann algebras},  International Congress of Mathematicians. 
Vol. I,  445--477, 
Eur. Math. Soc., Z$\ddot{\text{u}}$rich, 2007.
\item {[Po09]} S. Popa {\it Revisiting some problems in W*-rigidity}, available at 

http://www.math.ucla.edu/~popa/workshop0309/slidesPopa.pdf.
\item {[Wi08]} D. Witte Morris: {\it Introduction to Arithmetic Groups}, lecture notes, available at http://people.uleth.ca/~dave.morris/.

\enddocument